\documentclass[11pt]{article}
\usepackage{amsmath,amssymb}
\newcommand{\doublespace}
   {\addtolength{\baselineskip}{0.15\baselineskip}}

\newtheorem{pdef}{Definition}[section] %
\newtheorem{thm}[pdef]{Theorem}        

\newcounter{equationnumber}
\renewcommand{\theequation}{\thesection.\arabic{equation}}
\def\mathletters{
    \addtocounter{equation}{1}
    \edef\@currentlabel{\theequation}
    \setcounter{equationnumber}{\value{equation}}
    \setcounter{equation}{0}
    \edef\theequation{\@currentlabel\noexpand\alph{equation}}
    }

\title{A slight improvement to Korenblum's constant}
\author{Chun-Yen Shen}

\date{ }
\begin{document}
\maketitle \doublespace
\pagestyle{myheadings} \thispagestyle{plain} \markboth{   }{ }
\begin{abstract}
Let $A^2(D)$ be the Bergman space over the open unit disk $D$ in
the complex plane. Korenblum conjectured that there is an absolute
constant $c \in (0,1)$ such that whenever $|f(z)|\le |g(z)|$ in
the annulus $c<|z|<1$ then $||f(z)|| \le ||g(z)||$. This
conjecture had been solved by Hayman $[H1]$, but the constant $c$
in that paper is not optimal. Since then, there are many papers
dealing with improving the upper and lower bounds for the best
constant $c$. For example, in 2004 C.Wang gave an upper bound on
$c$,that is, $c < 0.67795$, and in 2006 A.Schuster gave a lower
bound ,$c > 0.21 $ .In this paper we slightly improve the upper
bound for $c$.
\end{abstract}
\footnotetext[1]{{\it 2000 Mathematics Subject Classification:}\,
Primary 30C80.} \footnotetext[2]{{\it Key words and
phrases.}\, Bergman space, Korenblum's constant.}
\section{Introduction}
Let D be the open unit disk
in the complex plane C.The Bergman space $A^2(D)$ consists of
analytic functions f in D such that

\[
         ||f||= \big[\int_{D} |f(z)|^2dA(z)\big]^{1/2} < \infty
      \]
where \[
      dA(z)={\frac{1}{\pi}}dxdy={\frac{1}{\pi}}rdrd\theta , z=x+iy=re^{i\theta}
       \]
is the normalized Lebesgue area measure on D. Korenblum
conjectured that there is an absolute constant c ,$0<c<1$,such
that whenever $|f(z)|\le|g(z)|$ in the annulus $c<|z|<1$ ,then
$||f(z)|| \le||g(z)||$. This conjecture is very natural and
inspired many work.The answer to this conjecture was obtained by Hayman.
He proved that the constant $c$ exists and is greater than $0.04$. But we can ask
whether this bound is optimal, and many papers have been done to
find the better upper and lower bounds for the constant $c$ (see
[H2],[S] and [W]). The best upper bound for $c$ until now is
$0.67795....$ in [W], and the best lower bound for this constant
is $0.21$ in [S]. In this paper we use Wang's example in a slightly more sophisticated way to improve the upper bound to
$0.677905$.
\section{Results}
The following theorem shows that $c < 0.677905 $.

\begin{thm}
Let
\[
f(z)=\frac{a+z^n}{2-az^n} , g(z)=\frac{z(1+az^n)}{ 2-az^n},
\]
where a=0.6666714 and n=10.Then $||f(z)||>||g(z)||$ and $|f(z)|\le
|g(z)|$ in $c<|z|<1$ ,where c=0.6779049274... is the real root of
the equation $f(z)=g(z)$
\end{thm}
\begin{pf} Let
\[
h(r)=max_{|z|=r}|\frac {f(z)}{g(z)}|=\frac{a+r^n}{r(1+ar^n)}.
\]
Then,$h(c)=h(1)=1.$ Since $f(z)/g(z)$ is analytic in $c\le |z| \le
1$,the maximum modulus theorem implies that $|f(z)|\le |g(z)|$ in
$c<|z|<1$. Using Maple to solve the equation ,we obtain that
c=0.6779049274....Next, a direct calculation shows that
\[
\int_{D} |f(z)|^2dA(z)-\int_{D} |g(z)|^2dA(z)
\]

\[
=\frac{1}{\pi}\int_0^{2\pi}\int_0^{1} \frac{a^2+2ar^ncosn\theta
+r^{2n}}{4-4ar^ncosn\theta +a^2r^{2n}}rdrd\theta -\frac{1}{\pi}\int_0^{2\pi}\int_0^{1}\frac{1+2ar^ncosn\theta +a^2r^{2n}}{4-4ar^ncosn\theta +a^2r^{2n}}r^3drd\theta
\]
Let $r^n=\rho$ and $n\theta =\phi$,then
\[
=\frac{1}{n\pi}\int_0^{2\pi}\int_0^{1}\frac{a^2+2a\rho cos\phi
+\rho^2}{4-4a\rho cos\phi+a^2\rho^2}\rho^{2/n-1}d\rho
d\phi-\frac{1}{n\pi}\int_0^{2\pi}\int_0^{1}\frac{1+2a\rho cos\phi
+a^2\rho^2}{4-4a\rho cos\phi+a^2\rho^2}\rho^{4/n-1}d\rho
d\phi
\]
Using Maple, we obtain that when a=0.6666714, n=10 and
c=0.6779049274...
Then
\[
\int_{D} |f(z)|^2dA(z)-\int_{D} |g(z)|^2dA(z)\ge 0.22*10^{-6}
\]
Hence
\[
||f(z)||>||g(z)||
\]
\end{pf}
\addcontentsline{toc}{section}{\protect \numberline {REFERENCES}}

\end{document}